# Separable equivalence, complexity and representation type

Simon F. Peacock

**Abstract.** We generalise the notion of separable equivalence, originally presented by Linckelmann in [Lin11b], to an equivalence relation on additive categories. We use this generalisation to show that from an initial equivalence between two algebras we may build equivalences between many related categories. We also show that separable equivalence preserves the representation type of an algebra. This generalises Linckelmann's result in [Lin11b], where he showed this in the case of symmetric algebras. We use these theorems to show that the group algebras of several small cyclic groups cannot be separably equivalent. This gives several examples of algebras that have the same complexity but are not separably equivalent.

[Lin11b] Linckelmann, *Finite generation of Hochschild cohomology of Hecke algebras of finite classical type in characteristic zero*, Bull. Lond. Math. Soc. **43** (2011), no. 5, 871–885

## 1  Introduction

In representation theory there are several notions of equivalence for algebras that can be used to help us understand certain algebras by comparing them to others. The most obvious equivalence we can describe is *isomorphism*, which from the point of view of representation theory is extremely rigid. More often we only care whether or not the module structure for the algebra is the same, which leads to the idea of *Morita equivalence*. We say that two algebras are Morita equivalent if their module categories are equivalent (as additive categories). From the module category for an algebra we can define further categories: the derived module category and the stable module category. Each of these give us a further equivalence relation on algebras, namely derived equivalence and stable equivalence, which are simply determined via (triangulated) equivalence between the relevant categories. Linckelmann presented the notion of *separable equivalence* in [Lin11b], which for self-injective algebras can be considered a generalisation of the other relations we have mentioned.

We begin the next section with Linckelmann's original definition of separable equivalence and provide some preliminary properties of the equivalence. We conclude the section by proving that complexity is preserved under separable equivalence

*Date*: January 11, 2017
*Key words and phrases.* separable equivalence; complexity; representation type
This research was supported by EPSRC funding

**Theorem 1:**

> Let $A$ and $B$ be finite dimensional algebras over a field $k$. If $A$ and $B$ are separably equivalent then their complexities agree:
>
> $$\operatorname{cx}(A) = \operatorname{cx}(B).$$

This property naturally leads us to to ask whether the converse is true: can we find algebras that have the same complexity but are not separably equivalent? Linckelmann gave the first example of such a situation in [Lin11a]. We will generalise this result in the final section where we show that group algebras of certain cyclic groups are not separably equivalent despite all having complexity 1. These results are corollaries to the main theorem of that section:

[Lin11a] Linckelmann, *Cohomology of block algebras of finite groups*, Representations of algebras and related topics, EMS Ser. Congr. Rep., Eur. Math. Soc., Zürich, 2011, pp. 189–250

**Theorem 8:**

> Let $\Lambda_n$ denote the truncated polynomial algebra $k[x]/(x^n)$ over an algebraically closed field $k$.
>
> The algebras $\Lambda_n$ and $\Lambda_m$ are not separably equivalent for positive integers $n \leq 6$ and $m \neq n$.

In order to prove this theorem we first require the results of sections 4 and 5. In section 4 we extend the definition of separable equivalence to allow for the notion of separable equivalence of categories. Under this extended definition we will see that two algebras are separably equivalent if and only if their module categories are separably equivalent. We can use this new definition to show that if we begin with a separable equivalence of algebras we may construct new categories that must also be separably equivalent. We additionally present a stronger form of separable equivalence, *symmetrical separably equivalence*, and show that for separably equivalent symmetric algebras, such as group algebras, we may always assume we have this stronger equivalence. These ideas provide the machinery required for the following two theorems:

**Theorem 2:**

> Let $A$ and $B$ be a pair of symmetrically separably equivalent $k$-algebras. Let $C$ and $D$ be a second pair of symmetrically separably equivalent $k$-algebras and let $\mathcal{E}$ be a small $k$-category. We have the following symmetrical separable equivalences of functor categories:
>
> (a) $\operatorname{Fun}(\operatorname{mod} A, \operatorname{mod} C) \sim \operatorname{Fun}(\operatorname{mod} B, \operatorname{mod} D)$,
>
> (b) $\operatorname{Fun}(\operatorname{mod} A, \mathcal{E}) \sim \operatorname{Fun}(\operatorname{mod} B, \mathcal{E})$,
>
> (c) $\operatorname{Fun}(\mathcal{E}, \operatorname{mod} A) \sim \operatorname{Fun}(\mathcal{E}, \operatorname{mod} B)$.

**Theorem 3:**

> Let $\mathcal{A}$ and $\mathcal{B}$ be symmetrically separably equivalent categories via $(F, G)$. If we have full subcategories $\mathcal{A}' < \mathcal{A}$ and $\mathcal{B}' < \mathcal{B}$ such that $F\mathcal{A}' \subseteq \mathcal{B}'$ and $G\mathcal{B}' \subseteq \mathcal{A}'$ then we have the following symmetrical separable equivalences:
>
> (a) $\mathcal{A}' \sim \mathcal{B}'$;
>
> (b) $\mathcal{A}/\mathcal{A}' \sim \mathcal{B}/\mathcal{B}'$;
>
> (c) $\operatorname{Fun}\left(\mathcal{A}/\mathcal{A}', \operatorname{mod} k\right) \sim \operatorname{Fun}\left(\mathcal{B}/\mathcal{B}', \operatorname{mod} k\right)$.

Then in section 5 we show that separable equivalence preserves the representation type of an algebra, which we use in the proof of theorem 8. In that proof we begin with a pair of algebras, which we wish to show are inequivalent and from these we build a pair of algebras that have different representation types. This then shows the starting pair cannot have been separably equivalent.

**Theorems 6 and 7:**

> Let $A$ and $B$ be finite dimensional algebras over an algebraically closed field $k$ such that $A$ and $B$ are separably equivalent.
>
> (a) If $A$ is of finite representation type then $B$ is of finite representation type.
>
> (b) If $A$ is a domestic algebra then $B$ is a domestic algebra.
>
> (c) If $A$ is an algebra of polynomial growth then $B$ is an algebra of polynomial growth.
>
> (d) If $A$ is of tame representation type then $B$ is of tame representation type.
>
> (e) If $A$ is of wild representation type then $B$ is of wild representation type.

**Acknowledgment.** I would like to thank Jeremy Rickard for all of his help over the past several years, without his guidance none of this work could exist.

## 2  Preliminaries

**Definition** (Separable equivalence [Lin11b])**.** Let $A$ and $B$ be finite dimensional algebras over a field $k$. We say that $A$ and $B$ are *separably equivalent* if there are bimodules ${}_A M_B$ and ${}_B N_A$ such that

(a) the modules ${}_A M$, $M_B$, ${}_B N$ and $N_A$ are finitely generated and projective; and

[Lin11b] Linckelmann, *Finite generation of Hochschild cohomology of Hecke algebras of finite classical type in characteristic zero*, Bull. Lond. Math. Soc. **43** (2011), no. 5, 871–885

(b) there are bimodules ${}_A X_A$ and ${}_B Y_B$ and bimodule isomorphisms

$$_A M \underset{B}{\otimes} N_A \xrightarrow{\sim} {}_A A_A \oplus {}_A X_A \qquad {}_B N \underset{A}{\otimes} M_B \xrightarrow{\sim} {}_B B_B \oplus {}_B Y_B$$

*Remarks.*

- If the bimodules $X$ and $Y$ are the zero modules then we would have a Morita equivalence.

- If the bimodules $X$ and $Y$ are projective (as bimodules) then we would have a stable equivalence (of Morita type).

- Separable equivalence can be considered a generalisation of the other equivalence relations as, for self-injective algebras,

$$\begin{aligned} \text{Isomorphism} &\implies \text{Morita equivalence} \\ &\implies \text{derived equivalence} \\ &\implies \text{stable equivalence of Morita-type} \\ &\implies \text{separable equivalence} \end{aligned}$$

The terminology *separable equivalence* comes from the following proposition that was stated by Linckelmann in [Lin11b].

**Proposition.** *A finite dimensional algebra $A$ over a field $k$ is separable (in the sense of [DI71]) if and only if it is separably equivalent to $k$.*

*Proof.* Firstly, assume that $A$ is a separable algebra so that (by definition) $A$ is a summand of $A \underset{k}{\otimes} A$ as $A$-$A$–bimodules. Taking $M = {}_A A_k$ and $N = {}_k A_A$, so that $M \underset{k}{\otimes} N = A \underset{k}{\otimes} A$, we have the required isomorphisms.

Now assume that $A$ is separably equivalent to $k$ through bimodules ${}_A M_k$ and ${}_k N_A$. Consider the functor $\operatorname{Hom}_{A-A}(M \underset{k}{\otimes} N, -)$:

$$\operatorname{Hom}_{A-A}\left(M \underset{k}{\otimes} N, -\right) \cong \operatorname{Hom}_A\left(M, \operatorname{Hom}_A(N, -)\right)$$
$$= \operatorname{Hom}_A(M, -) \circ \operatorname{Hom}_A(N, -)$$

Since $M$ is projective as a left $A$–module we have that $\operatorname{Hom}_A(M, -)$ is exact. Similarly the functor $\operatorname{Hom}_A(N, -)$ is exact and hence so is the composition. We therefore have that $M \underset{k}{\otimes} N$ is projective as an $A$-$A$–bimodule and so $A$ is projective as an $A$-$A$–bimodule.

□

[Lin11b] Linckelmann, *Finite generation of Hochschild cohomology of Hecke algebras of finite classical type in characteristic zero*, Bull. Lond. Math. Soc. **43** (2011), no. 5, 871–885

[DI71] DeMeyer and Ingraham, *Separable algebras over commutative rings*, Lecture Notes in Mathematics, vol. 181, Springer-Verlag, Berlin-New York, 1971

*Remark.* If we had assumed $A$ were an algebra over a ring $R$, rather than over a field, we would not have that $R$ is a summand of $A$ and so the above proof does not go through. If we had this additional assumption however, (for example if $A$ were free as an $R$-module) we could follow the same proof.

It will be convenient to talk about situations in which only one of the isomorphisms in the definition of separable equivalence exists. In this situation we will use the language of Bergh and Erdmann (see [BE11, 2506f.]) and say that one algebra *separably divides* the other.

**Definition** (Separably divides)**.** Given two $R$-algebras $A$ and $B$, we say that *A separably divides B* if there exist bimodules ${}_AM_B$ and ${}_BN_A$, finitely generated projective on both sides, such that $A$ is a bimodule direct summand of $M \underset{B}{\otimes} N$.

*Remark.* Notice that two algebras $A$ and $B$ are separably equivalent if and only if $A$ separable divides $B$ and $B$ separably divides $A$.

The next proposition, which was stated by Linckelmann in [Lin11b], will give us our first example of algebras that are separably equivalent but are not equivalent in any of the more specialised ways we have mentioned.

**Proposition.** *Let $G$ be a finite group, $k$ a field of characteristic $p > 0$. If $P$ is a Sylow $p$-subgroup of $G$ then $kP$ is separably equivalent to $kG$.*

*Proof.* The separable equivalence is given by the bimodules ${}_{kP}kG_{kG}$ and ${}_{kG}kG_{kP}$ and the proof follows an argument similar to the proof of Mackey's decomposition theorem. □

*Example.* Let $A_4$ denote the alternating group. The Sylow 2-subgroups of $A_4$ are isomorphic to the Klein 4-group, $V_4$ and so the proposition above tells us that over a field $k$ of characteristic 2, $kA_4$ and $kV_4$ are separably equivalent. Similarly the Sylow 2-subgroups of $A_5$ are isomorphic to $V_4$ and so $kA_5$ is also separably equivalent to $kV_4$, and hence all three algebras are separably equivalent. An application of Green correspondence demonstrates that $kA_4$ and $kA_5$ are stably equivalent (see [Alp86, theorem 10.1] for example) however neither algebra is stably equivalent to $kV_4$. This demonstrates that separable equivalence is not a reformulation of any of the other equivalences and therefore gives us something new to work with.

In the proof of the above proposition the bimodules that we used to form the separable equivalence were duals of one another. In particular we have that tensoring with $M$ is both left and right adjoint to tensoring with $N$. If $A$ and $B$ are symmetric separably equivalent algebras then we may always choose the modules so that these adjunctions exist. We have the following definition and proposition.

**Definition** (Symmetrical separable equivalence)**.** Let $A$ and $B$ be finite dimensional algebras. We say that $A$ and $B$ are *symmetrically separably equivalent* if there is a separable equivalence $(M, N)$ such that $- \otimes M$ is both left and right adjoint to $- \otimes N$.

[BE11] Bergh and Erdmann, *The representation dimension of Hecke algebras and symmetric groups*, Adv. Math. **228** (2011), no. 4, 2503–2521

[Lin11b] Linckelmann, *Finite generation of Hochschild cohomology of Hecke algebras of finite classical type in characteristic zero*, Bull. Lond. Math. Soc. **43** (2011), no. 5, 871–885

[Alp86] Alperin, *Local representation theory*, Cambridge Studies in Advanced Mathematics, vol. 11, Cambridge University Press, Cambridge, 1986

**Proposition.** *If A and B are symmetric algebras then A and B are separably equivalent if and only if A and B are symmetrically separably equivalent.*

*Proof.* If $M, N$ form a separable equivalence then $M \oplus N^*$ and its dual form a symmetrical separable equivalence. □

*Remark.* The above proposition shows that for symmetric algebras the two definitions, separable equivalence and symmetrical separable equivalence, coincide. Whether or not these definitions coincide in the general case is an open question.

## 3   Complexity

There are certain properties of algebras that are preserved through separable equivalence. One such property is the complexity of an algebra. Here we introduce what is meant by the complexity of a module and of an algebra and go on to prove that it is unchanged by separable equivalence. This result seems to be well-known but we don't know of any complete statement or proof in print.

**Definition** (Complexity). Let $A$ be an algebra over a field $k$, let $M$ be an $A$-module and
$$\cdots \longrightarrow P_1 \longrightarrow P_0 \longrightarrow M \longrightarrow 0$$
a projective resolution of $M$, which we will denote by $P_*$.

If there exists an integer $d$, such that for some $\lambda \in \mathbb{N}$ we have $\dim(P_n) \leq \lambda n^{d-1}$ for all $n \in \mathbb{N}$ then we say that $P_*$ has *finite complexity* and we call the smallest such $d$ the *complexity of the resolution*, which we denote by $\operatorname{cx} P_*$.

The *complexity of the module $M$* is equal to the complexity of a minimal projective resolution of $M$.

The *complexity of an algebra* is the maximal complexity for a module of that algebra.

*Remark.* If $P_* \to M$ is a minimal projective resolution of $M$ and $Q_* \to M$ is any other projective resolution then $\operatorname{cx}(P_*) \leq \operatorname{cx}(Q_*)$. In particular $\operatorname{cx}(M) = \operatorname{cx}(P_*)$.

We require two lemmas before we can prove the main result of this section: that separable equivalence preserves complexity.

**Lemma 3.1.** *Let A and B be finite dimensional algebras over a field k and $_A M_B$ a bimodule that is finitely generated projective as a both an A-module and as a B-module. If we have a projective resolution of A-modules*
$$\cdots \longrightarrow P_1 \longrightarrow P_0 \longrightarrow X \longrightarrow 0$$
*then $\operatorname{cx}(X \underset{A}{\otimes} M) \leq \operatorname{cx}(X) \leq \operatorname{cx}(P_*)$.*

*Proof.* Since $M$ is projective as a $A$-module the functor

$$- \underset{A}{\otimes} M \colon \operatorname{mod} A \longrightarrow \operatorname{mod} B$$

is exact and since $M$ is projective as a $B$-module each $P_i \underset{A}{\otimes} M$ is projective. We therefore have that $P_* \underset{A}{\otimes} M$ is a projective resolution of $X \underset{A}{\otimes} M$.

It is clear that $\dim(P_i \underset{A}{\otimes} M) \leq \dim(P_i \underset{k}{\otimes} M)$ and since $M$ and $B$ are both finitely generated $\dim(P_i \underset{k}{\otimes} M) = \dim(P_i)\dim(M) < \infty$.

Finally if $\dim(P_i) \leq f(i)$ for some polynomial $f$ then $\dim(P_i \underset{A}{\otimes} M) \leq \dim(M)f(i)$ and hence $\operatorname{cx}(P_* \underset{A}{\otimes} M) \leq \operatorname{cx}(P_*)$. □

**Lemma 3.2.** *If $A$ separably divides $B$ via the modules $({}_A M_B, {}_B N_A)$ and $X_A$ is an $A$-module then $\operatorname{cx}(X \underset{A}{\otimes} M) = \operatorname{cx}(X)$.*

*Proof.* It suffices to show that $\operatorname{cx}(X) \leq \operatorname{cx}(X \underset{A}{\otimes} M \underset{B}{\otimes} N)$ as together with lemma 3.1 this gives

$$\operatorname{cx}(X \underset{A}{\otimes} M \underset{B}{\otimes} N) \leq \operatorname{cx}(X \underset{A}{\otimes} M) \leq \operatorname{cx}(X) \leq \operatorname{cx}(X \underset{A}{\otimes} M \underset{B}{\otimes} N)$$

and we will have equality throughout. Since $X \underset{A}{\otimes} M \underset{B}{\otimes} N \cong X \oplus X'$ for some $X'$ we have that a minimal projective resolution of $X \underset{A}{\otimes} M \underset{B}{\otimes} N$ is simply the direct sum of the minimal projective resolutions of $X$ and $X'$ and thus the given inequality is immediate. □

The next theorem is a direct consequence of the previous three lemmas.

**Theorem 1:**

> If $A$ separably divides $B$ then $\operatorname{cx}(A) \leq \operatorname{cx}(B)$.
> 
> If $A$ and $B$ are separably equivalent then $\operatorname{cx}(A) = \operatorname{cx}(B)$.

If $G$ is a finite group and $k$ a field of characteristic $p$ then results of Alperin and Evens in [AE81] show that the complexity of $kG$ is equal to the $p$-rank of $G$, that is the rank of the largest elementary abelian $p$-subgroup of $G$. For example if $G$ is a cyclic $p$-group then $\operatorname{cx} kG = 1$.

[AE81] Alperin and Evens, *Representations, resolutions and Quillen's dimension theorem*, J. Pure Appl. Algebra **22** (1981), no. 1, 1–9

The remarks above lead us naturally to ask if all group algebras of cyclic $p$-groups are separably equivalent. We will answer this question in the negative in section 6, but in order to achieve the results therein we must first broaden our definition of separable equivalence to abstract categories. If we begin with a pair of separably equivalent algebras we can use this broadened definition to build additional pairs of separably equivalent algebras (or categories) from them.

# 4 Categorical formulation

The notion of separable equivalence we have used thus far was presented in terms of bimodules for two algebras but note that we could have just as easily defined this equivalence in terms of functors between the module categories. If $A$ and $B$ are separably equivalent algebras via $(M, N)$ then we may define functors $F = - \underset{A}{\otimes} M$ and $G = - \underset{B}{\otimes} N$. Notice that these functors are exact, send projective modules to projective modules, and their composition (in either order) contains the relevant identity functor as a summand. In fact given functors $F \colon \operatorname{mod} A \to \operatorname{mod} B$ and $G \colon \operatorname{mod} B \to \operatorname{mod} A$ with these properties we may define the modules $M = FA$ and $N = GB$ and this pair gives rise to a separable equivalence. Together, these remarks allow us to generalise the definition of separable equivalence to exact categories.

**Definition** (Separable equivalence). Let $\mathcal{A}$ and $\mathcal{B}$ be exact categories. We say that these categories are separably equivalent if there are exact (additive) functors such that

$$\mathcal{A} \underset{G}{\overset{F}{\rightleftarrows}} \mathcal{B}$$

- the image of a projective object is projective;
- the identity functor is a summand of $GF$ and of $FG$.

In a similar way we may generalise the definition of symmetrical separable equivalence and in this case we may even drop the requirement that the categories are exact.

**Definition** (Symmetrical separable equivalence). Let $\mathcal{A}$ and $\mathcal{B}$ be additive categories. We say that these categories are *symmetrically separably equivalent* if there are additive functors such that

$$\mathcal{A} \underset{G}{\overset{F}{\rightleftarrows}} \mathcal{B}$$

- both $(F, G)$ and $(G, F)$ form an adjunction; and
- the identity functor is a summand of $GF$ and of $FG$.

*Remark.* If $\mathcal{A}$ and $\mathcal{B}$ are module categories then the adjointness implies that the functors are exact and that projective modules are sent to projective modules.

With these generalised definitions we are in a position to demonstrate how we may begin with a separable equivalence of algebras and build further equivalences from this starting point. We begin with a proposition using the original bimodule definition.

**Proposition.** *Let A, B and C be algebras over a field k. If A separably divides B then $A \otimes_k C$ separably divides $B \otimes_k C$.*

*Proof.* Let $_A M_B$ and $_B N_A$ be a pair of bimodules with the property that

$$_A M \otimes_B N_A \cong {}_A A_A \oplus {}_A X_A.$$

The tensor product $M \otimes_k C$ is an $(A \otimes_k C)$-$(B \otimes_k C)$–bimodule with the actions

$$(a \otimes c_1)(m \otimes c_2)(b \otimes c_3) = amb \otimes c_1 c_2 c_3$$

and we can similarly define a $(B \otimes_k C)$-$(A \otimes_k C)$ action on $N \otimes_k C$.

Thus we have that

$$(M \otimes_k C) \otimes_{B \otimes C} (N \otimes_k C) \cong (M \otimes_B N) \otimes_k C$$
$$\cong (A \oplus X) \otimes_k C$$
$$\cong (A \otimes_k C) \oplus (X \otimes_k C)$$

Now if $_A M$ is projective then $_A M$ is a summand of $A^n$ for some $n$. Therefore $M \otimes_k C$ is a summand of

$$A^n \otimes_k C \cong (A \otimes_k C)^n$$

and hence is projective. □

In the remainder of this section we will limit ourselves to the case where we have a symmetrical separably equivalence, and so we begin with a lemma regarding adjoint functors.

**Lemma 4.1.** *Let $\mathcal{A}$, $\mathcal{B}$, $\mathcal{C}$ and $\mathcal{D}$ be k-categories. Given functors*

$$\mathcal{A} \underset{R}{\overset{L}{\rightleftarrows}} \mathcal{B} \qquad \mathcal{C} \underset{G}{\overset{F}{\rightleftarrows}} \mathcal{D}$$

*we can define the functors*

$$\mathrm{Fun}(\mathcal{A},\mathcal{C}) \underset{G \circ \_ \circ L}{\overset{F \circ \_ \circ R}{\rightleftarrows}} \mathrm{Fun}(\mathcal{B},\mathcal{D})$$

*where* $\mathrm{Fun}(\mathcal{A},\mathcal{C})$ *denotes the category of functors* $\mathcal{A} \to \mathcal{C}$.

*If $(L, R)$ and $(F, G)$ are both adjoint pairs then so is $(F\text{—}R, G\text{—}L)$.*

*Proof.* Let $\operatorname{Nat}(S, T)$ denote the category of natural transformations $S \Rightarrow T$.

The adjunction isomorphisms are given by

$$\operatorname{Nat}(FXR, Y) \xrightarrow{\Phi_{X,Y}} \operatorname{Nat}(X, GYL)$$
$$\xi \mapsto G\xi L \circ \eta^{FG} X \eta^{LR}$$

and

$$\operatorname{Nat}(X, GYL) \xrightarrow{\Psi_{X,Y}} \operatorname{Nat}(FXR, Y)$$
$$\nu \mapsto \varepsilon^{FG} Y \varepsilon^{LR} \circ F\nu R$$

where $\eta$ and $\varepsilon$ are the units and counits of the adjunctions with superscripts indicating the adjunction in question.

One direction of the proof that these are inverse mappings is demonstrated in the following commutative diagram.

$$\begin{array}{c}
\xymatrix{
& & FXR \ar[rr]^{\xi} & & Y \\
FXR \ar[urr]^{\mathrm{Id}} \ar[dr]_{F\eta^{FG}XR} \ar[dd]_{F\eta^{FG}XR} & & & & \\
& & FGFXR \ar[u]^{\varepsilon^{FG}FXR} \ar[rr]_{FG\xi} & & FGY \ar[u]_{\varepsilon^{FG}Y} \\
FGFXR \ar[urr]^{\mathrm{Id}} \ar[dr]_{FGFX\eta^{LR}R} & & & & \\
& & FGFXRLR \ar[u]^{FGFXR\varepsilon^{LR}} \ar[rr]_{FG\xi LR} & & FGYLR \ar[u]_{FGY\varepsilon^{LR}}
}
\end{array}$$

The composition $\Psi\Phi$ is given by the path from $FXR$ to $Y$ around the exterior of the diagram. The squares commute because of the natural transformations and the triangles commute because of the adjunctions. □

As a direct result of this lemma we have the following theorem.

**Theorem 2:**

> Let $A$ and $B$ be a pair of symmetrically separably equivalent $k$-algebras. Let $C$ and $D$ be a second pair of symmetrically separably equivalent $k$-algebras and let $\mathcal{E}$ be a small $k$-category. We have the following symmetrical separable equivalences of functor categories:

  (a) $\mathrm{Fun}(\mathrm{mod}\, A, \mathrm{mod}\, C) \sim \mathrm{Fun}(\mathrm{mod}\, B, \mathrm{mod}\, D)$,

  (b) $\mathrm{Fun}(\mathrm{mod}\, A, \mathcal{E}) \sim \mathrm{Fun}(\mathrm{mod}\, B, \mathcal{E})$,

  (c) $\mathrm{Fun}(\mathcal{E}, \mathrm{mod}\, A) \sim \mathrm{Fun}(\mathcal{E}, \mathrm{mod}\, B)$.

If we consider the restriction of the functors that make up a separable equivalence we may find that some subcategories are also equivalent. For instance if we have a symmetrical separable equivalence $(F, G)$ between two categories $\mathcal{A}$ and $\mathcal{B}$, such that $(F, G)$ restrict to functors between full subcategories $\mathcal{A}'$ and $\mathcal{B}'$ then in fact we have a symmetrical separable equivalence of these subcategories. This follows simply from the fact that the identity functor on the subcategory is the restriction of that of the parent category and additionally that, for full subcategories, the Hom-sets are equal to those of the parent category. This together with a dual result gives us:

**Theorem 3:**

Let $\mathcal{A}$ and $\mathcal{B}$ be symmetrically separably equivalent categories via $(F, G)$. If we have full subcategories $\mathcal{A}' < \mathcal{A}$ and $\mathcal{B}' < \mathcal{B}$ such that $F\mathcal{A}' \subseteq \mathcal{B}'$ and $G\mathcal{B}' \subseteq \mathcal{A}'$ then we have the following symmetrical separable equivalences:

  (a) $\mathcal{A}' \sim \mathcal{B}'$;

  (b) $\mathcal{A}/\mathcal{A}' \sim \mathcal{B}/\mathcal{B}'$;

  (c) $\mathrm{Fun}\left(\mathcal{A}/\mathcal{A}', \mathrm{mod}\, k\right) \sim \mathrm{Fun}\left(\mathcal{B}/\mathcal{B}', \mathrm{mod}\, k\right)$.

*Remark.* Here $\mathcal{A}/\mathcal{A}'$ denotes the category whose objects are the objects of $\mathcal{A}$ and morphisms are given by equivalence classes of morphisms in $\mathcal{A}$ under the relation: $f, g \in \mathrm{Hom}_\mathcal{A}(X, Y)$ then $f \sim g$ if and only if $f - g$ factors through an object of $Z \in \mathcal{A}'$.

$$X \xrightarrow{f-g} Y \searrow \nearrow Z$$

*Proof.* The proof of (a) is clear from the discussion above and (b) from a dual argument. For (c) we need only note that

$$\mathrm{Fun}\left(\mathcal{A}/\mathcal{A}', \mathrm{mod}\, k\right)$$

is a full subcategory of $\mathrm{Fun}(\mathcal{A}, \mathrm{mod}\, k)$ via the embedding that composes a functor with the obvious projection

$$\mathcal{A} \longrightarrow \mathcal{A}/\mathcal{A}'$$

and that the equivalence given by theorem 2 restricts to functors of these subcategories.

□

*Example.* If $A$ and $B$ are separably equivalent symmetric algebras then the subcategories of projective modules, proj $A$ and proj $B$, are symmetrically separably equivalent. Similarly the category of representations of their stable categories, $\text{Fun}(\underline{\text{mod}}\, A, \text{mod}\, k)$ and $\text{Fun}(\underline{\text{mod}}\, B, \text{mod}\, k)$, are symmetrically separably equivalent.

## 5   Representation type

Recall that each algebra has a property referred to as its representation type. For an algebra over an algebraically closed field Drozd showed in [Dro77] and [Dro80] that there are only three possibilities for its representation type: finite, tame or wild. The purpose of this section is to prove that this property is preserved under separable equivalence. Note that this is a generalisation of a result by Linckelmann ([Lin11b, proposition 3.3]), which proved this in the case of symmetric algebras. In section 6 we will use the results of this section to show that certain group algebras cannot be separably equivalent. At first it may seem that Linckelmann's result would be sufficient for our purpose, since group algebras are symmetric, however we will actually be showing the representation type differs for some of the constructions at the end of the last section. These algebras cannot be guaranteed to be symmetric, even if the original algebras were.

We begin with the definitions of finite, tame and wild representation types.

**Definition** (Finite representation type). An algebra $A$ is said to have *finite representation type* if there exists only finitely many isomorphism classes of indecomposable right (equivalently left) $A$-modules.

**Definition** (Tame representation type). An algebra $A$ over a field $k$ is said to have *tame representation type* if it does not have finite representation type and given any $d \in \mathbb{N}$ there is a finite set of $k[t]$-$A$–bimodules $\{X_i\}$, free and finitely-generated as $k[t]$-modules, such that for all but finitely many $d$-dimensional indecomposable $A$-modules $M$ (up to isomorphism) we have

$$M \cong \frac{k[t]}{(t - \lambda)} \underset{k[t]}{\otimes} X_i$$

for some $X_i$ and some $\lambda \in k$.

*Remarks.*

- We may interpret this definition as saying that an algebra is tame if its isomorphism classes of indecomposable $d$-dimensional modules can be classified by a finite number of 1-parameter families of modules, for each $d$.

[Dro77] Drozd, *Tame and wild matrix problems*, Matrix problems (Russian), Akad. Nauk Ukrain. SSR Inst. Mat., Kiev, 1977, pp. 104–114

[Dro80] Drozd, *Tame and wild matrix problems*, Representation theory, II (Proc. Second Internat. Conf., Carleton Univ., Ottawa, Ont., 1979), Lecture Notes in Math., vol. 832, Springer, Berlin, 1980, pp. 242–258

[Lin11b] Linckelmann, *Finite generation of Hochschild cohomology of Hecke algebras of finite classical type in characteristic zero*, Bull. Lond. Math. Soc. **43** (2011), no. 5, 871–885

- Some authors include finite representation type within the class of tame representation type however here we consider the two to be mutually exclusive.

**Definition** (Wild representation type). An algebra $A$ over a field $k$ is said to have *wild representation type* if there is a $k\langle u, v\rangle$-$A$–bimodule $X$, finitely generated free as a $k\langle u, v\rangle$-module, such that

- if $M$ is an indecomposable right $k\langle u, v\rangle$-module then $M \underset{k\langle u,v\rangle}{\otimes} X$ is indecomposable;
- for $k\langle u, v\rangle$-modules $M$ and $N$: if $M \underset{k\langle u,v\rangle}{\otimes} X \cong N \underset{k\langle u,v\rangle}{\otimes} X$ then $M \cong N$.

Here the notation $k\langle u, v\rangle$ represents the free $k$-algebra on two generators.

*Remark.* An algebra has wild representation type if its module category is at least as complicated as the module category for the free algebra on two generators.

The representation type of an algebra is intimately linked to what are known as *generic modules*. In fact the three definitions of representation type can be restated purely in terms of these modules, and we will use these definitions to prove the preservation of representation type.

Recall that if $M_A$ is a right $A$-module then $M$ is naturally a left module for its endomorphism ring $\text{End}(M_A)$ (in fact it is an $\text{End}(M_A)$-$A$–bimodule).

**Definition** (Endolength). Let $M_A$ be an $A$-module. We say that the *endolength* of $M$ is its length when considered as a module for its endomorphism ring and denote this by $\text{endlen}(M)$. We say that the module is *endofinite* if it has finite endolength.

**Definition** (Generic module). An indecomposable $A$-module $M$ is said to be a *generic module* if it has infinite length over $A$ but has finite length over $\text{End}(M_A)$.

Suppose $A$ is a tame algebra so that for each dimension $d \in \mathbb{N}$ we have a finite collection of $k[t]$-$A$–bimodules satisfying certain properties. We denote by

$$\mu_A(d) \in \mathbb{N}$$

the minimum number of these modules required to satisfy the definition. The following theorem gives the link between generic modules and representation type.

**Theorem 4:** *Crawley-Boevey,* [CB91, 5.7]

> For an algebra $A$ over an algebraically closed field, let $g_A(n)$ denote the number of isomorphism classes of generic $A$-modules of endolength $n$. Then
> 
> $$\mu_A(n) = \sum_{d|n} g_A(d)$$

[CB91] Crawley-Boevey, *Tame algebras and generic modules*, Proc. London Math. Soc. (3) **63** (1991), no. 2, 241–265

As immediate corollaries to theorem 4 we can provide alternative definitions of finite, tame and wild type in terms of the number of generic modules.

**Corollary.** *An algebra A is of finite type if and only if $g_A(n) = 0$ for all $n \in \mathbb{N}$.*

**Corollary.** *An algebra A over an algebraically closed field is tame if and only if $g_A(n) < \infty$ for all $n \in \mathbb{N}$ and $g_A(n) > 0$ for some $n \in \mathbb{N}$.*

**Corollary.** *An algebra A over an algebraically closed field is wild if and only if $g_A(n) = \infty$ for some $n \in \mathbb{N}$.*

In order to prove that separable equivalence preserves representation type we will need the following theorem on decomposability of endofinite modules.

**Theorem 5:** *Endofinite Decomposability*

> If $M$ is an endofinite module then there is a finite set of indecomposable modules $M_i$, and cardinals $\kappa_i$, such that $M$ decomposes as
> 
> $$M \cong M_1^{(\kappa_1)} \oplus \cdots \oplus M_n^{(\kappa_n)}.$$
> 
> Moreover if $M_i \not\cong M_j$ for all $i \neq j$ then $\operatorname{end \, len}(M) = \sum_{i=1}^{n} \operatorname{end \, len}(M_i)$.
> 
> The notation $M^{(\kappa)}$ denotes the direct sum of $\kappa$ copies of $M$.

For details on the proof of this theorem see [Pre09, 4.4.29], which is also a very good reference on endofinite modules in general.

We require one final lemma before we can show that representation type is preserved under separable equivalence.

**Lemma 5.1.** *Let ${}_A M_B$ be a finitely generated bimodule. There is a constant $c_M$ such that if $X_A$ is endofinite then $X \underset{A}{\otimes} M$ has endolength*

$$\operatorname{end \, len}\left(X \underset{A}{\otimes} M\right) \leq c_M \operatorname{end \, len}(X).$$

*Proof.* As $M$ is finitely generated there is an integer $n$ and a left $A$-epimorphism $A^n \to M$. This map gives an exact sequence of $\operatorname{End}_A(X)$ modules

$$\begin{array}{ccc} X \underset{A}{\otimes} A^n & \longrightarrow X \underset{A}{\otimes} M \longrightarrow 0 \\ \shortparallel & \\ X^n & \end{array}$$

Thus $n \operatorname{end \, len}(X)$ bounds the length of $X \underset{A}{\otimes} M$ as an $\operatorname{End}_A(X)$-module.

[Pre09] Prest, *Purity, spectra and localisation*, Encyclopedia of Mathematics and its Applications, vol. 121, Cambridge University Press, Cambridge, 2009

If we have a chain of $\mathrm{End}_B(X \underset{A}{\otimes} M)$-modules

$$0 = M_r \subsetneq \cdots \subsetneq M_1 \subsetneq X \underset{A}{\otimes} M$$

then this can be considered as a chain of $\mathrm{End}_A(X)$-modules via the canonical homomorphism

$$\begin{aligned} \mathrm{End}_A(X) &\longrightarrow \mathrm{End}_B(X \underset{A}{\otimes} M) \\ \phi &\mapsto \phi \underset{A}{\otimes} M \end{aligned}$$

and hence $\mathrm{end\,len}(X \underset{A}{\otimes} M) \leq n\, \mathrm{end\,len}(X)$. □

**Theorem 6:**

> Let $A$ and $B$ be finite dimensional algebras over an algebraically closed field $k$ such that $A$ separably divides $B$.
>
> (a) If $B$ is of finite representation type then $A$ is of finite representation type.
>
> (b) If $B$ is of tame representation type then $A$ is of finite or tame representation type.
>
> In particular the representation type of an algebra is preserved by separable equivalence.

*Proof.* We begin by proving part (b).

Let $_A M_B$ and $_B N_A$ be the modules providing the separable division.

Denote the generic $B$-modules of endolength $d$ by

$$^d G_1, {}^d G_2, \ldots, {}^d G_{g_B(d)}.$$

If $H$ is a generic $A$-module of endolength $d$ then by theorem 5 and lemma 5.1

$$H \underset{A}{\otimes} M \cong \bigoplus_{j=1}^{m} {}^{d_j} G_{i_j}^{(\kappa_j)} \oplus F$$

with $1 \leq m \leq c_M d$ and $d_j \leq c_M d$ for all $j$ and for some finite length module $F$. That this decomposition is essentially unique follows from section 4 of [CB92], in particular see the remarks following proposition 4.5.

We have that $H$ is a summand of $H \underset{A}{\otimes} M \underset{B}{\otimes} N$. If $H$ is a summand of $F \underset{B}{\otimes} N$ then $H$ has finite length, which is a contradiction as $H$ is generic. Therefore $H$ is a summand of $^{d_j} G_{i_j} \underset{B}{\otimes} N$ for some $j$.

[CB92] Crawley-Boevey, *Modules of finite length over their endomorphism rings*, Representations of algebras and related topics (Kyoto, 1990), London Math. Soc. Lecture Note Ser., vol. 168, Cambridge Univ. Press, Cambridge, 1992, pp. 127–184

Now define $\mathcal{H}(d)$ as follows:

$$\mathcal{H}(d) = \left\{ H \in \mathrm{mod}\, A \,\middle|\, \begin{array}{l} H \text{ a generic summand of } {}^{d'}G_i \underset{B}{\otimes} N \\ d' \leq c_M d,\ 1 \leq i \leq g_B(d') \end{array} \right\}.$$

Thus if $H$ is any generic module with $\mathrm{end\,len}(H) \leq d$ then $H \in \mathcal{H}(d)$.

We have assumed that $B$ is tame and therefore $g_B(n) < \infty$ for all $n \in \mathbb{N}$, in particular the set $\mathcal{H}(d)$ is finite and so $A$ must be of finite or tame representation type.

For part (a) there are no generic $B$-modules and so if $H$ is a generic $A$-module then $H \underset{A}{\otimes} M \cong F$ for some finite length module $F$. This would mean that $H$ is of finite length, a contradiction. □

Suppose $A$ is a tame algebra and we again use $\mu_A(d)$ to denote the minimum number of 1-parameter families of indecomposable $A$-modules. By imposing upper bounds on $\mu_A$ we can further subdivide the tame algebras. These subdivisions are called domestic algebras and algebras of polynomial growth.

**Definition** (Domestic algebra). An algebra $A$ is said to be *domestic* if there is some integer $N$ such that $\mu_A(n) \leq N$ for all positive integers $n$.

**Definition** (Polynomial growth). An algebra $A$ is said to be of *polynomial growth* if there are some positive integers $C$ and $\gamma_A$ such that

$$\mu_A(n) \leq C n^{\gamma_A}$$

for all positive integers $n$. If the integer $\gamma_A$ is chosen minimally with respect to the definition then it is called the *growth rate*.

We may again characterise these classes of algebras using the number of generic modules. As further corollaries to theorem 4 we have:

**Corollary.** *An algebra $A$ is of polynomial growth if and only if there are integers $C$ and $\delta$ such that*

$$g_A(n) \leq C n^\delta$$

*for all positive integers $n$.*

**Corollary.** *An algebra $A$ is domestic if and only if it has only finitely many generic modules.*

Now we can see that separable equivalence also preserves these subdivisions of tameness.

**Theorem 7:**

Let $A$ and $B$ be finite dimensional algebras over an algebraically closed field $k$ such that $A$ separably divides $B$.

(a) If $B$ is domestic then $A$ is domestic.

(b) If $B$ is of polynomial growth then $A$ is of polynomial growth.

In particular the properties of domestic and polynomial growth are preserved under separable equivalence.

*Proof.* Define $\mathcal{H}(d)$ as in the proof of theorem 6. For part (b) it is enough to bound the cardinality of $\mathcal{H}(d)$ by a polynomial in $d$.

The number of distinct endofinite summands of ${}^{d'}G_i \otimes_B N$ is bounded above by its endolength and hence by $c_N d'$. The number of generic modules of endolength $d'$ is given by $g_B(d')$. We have

$$|\mathcal{H}(d)| \leq \sum_{d' \leq c_M d} (c_N d') g_B(d')$$

$$\leq \sum_{d' \leq c_M d} (c_N d') C d'^{\delta}$$

$$\leq (c_M d)(c_N c_M d) C (c_M d)^{\delta}$$

$$\leq C' d^{\delta+2}$$

and hence $A$ is of polynomial growth.

For part (a) we must show that the number of generic $A$ modules is finite. As $B$ is domestic there are only finitely many generic modules, thus there is some integer $d$ for which $g_B(d') = 0$ for all $d' > d$. In particular $\mathcal{H}(d) = \mathcal{H}(d')$ for all $d' > d$. Thus every generic $A$ module is in the finite set $\mathcal{H}(d)$. □

## 6   Examples of inequivalence

The remarks at the end of section 3 led us to ask whether or not there exist cyclic groups $C_{p^n}$ and $C_{p^m}$ such that their group algebras over a field of characteristic $p$ are separably equivalent. Unfortunately the best we can offer is a partial solution to this problem. Specifically we will demonstrate the inequivalence of the group algebras for several small cyclic groups, leaving the general question wide open.

Let

$$\Lambda_n = \frac{k[x]}{(x^n)}$$

denote the truncated polynomial algebra of length $n$. Over a field of characteristic $p$ we have an isomorphism $kC_{p^n} \cong \Lambda_{p^n}$. By phrasing the above question in terms of $\Lambda_n$ we can achieve results for a general field that will then give answers for the cyclic group case over fields with the correct characteristic.

The Auslander–Reiten quiver for $\Lambda_n$ is

$$1 \underset{\beta_1}{\overset{\alpha_1}{\rightleftarrows}} 2 \underset{\beta_2}{\overset{\alpha_2}{\rightleftarrows}} \cdots \underset{\beta_{n-1}}{\overset{\alpha_{n-1}}{\rightleftarrows}} n \qquad \begin{array}{l} \alpha_1\beta_1 = 0 \\ \alpha_i\beta_i = \beta_{i-1}\alpha_{i-1} \text{ for } 1 < i < n \end{array}$$

If $kQ_n$ is the Auslander algebra of $\Lambda_n$ then the category $\operatorname{mod} kQ_n$ is equivalent to the category $\operatorname{Fun}(\operatorname{mod}\Lambda_n, \operatorname{mod} k)$. This fact together with theorem 2 means that if $\Lambda_n$ and $\Lambda_m$ are separably equivalent then $kQ_n$ and $kQ_m$ are separably equivalent.

If we instead restrict to the stable category, we must add additional relations for all maps that factor through a projective module. In this example we have a single projective module $\Lambda_n$, represented by vertex $n$ in the Auslander-Reiten quiver, thus we must add the relation $\beta_{n-2}\alpha_{n-2} = 0$ and remove the vertex $n$.

$$1 \underset{\beta_1}{\overset{\alpha_1}{\rightleftarrows}} 2 \underset{\beta_2}{\overset{\alpha_2}{\rightleftarrows}} \cdots \underset{\beta_{n-2}}{\overset{\alpha_{n-2}}{\rightleftarrows}} n-1 \qquad \begin{array}{l} \alpha_1\beta_1 = \beta_{n-2}\alpha_{n-2} = 0 \\ \alpha_i\beta_i = \beta_{i-1}\alpha_{i-1} \text{ for } 1 < i < n-1 \end{array}$$

The path algebra for the quiver with relations given above is called the preprojective algebra of type $A_{n-1}$. Here the $A_{n-1}$ refers to the Dynkin diagram and note that preprojective algebras can be defined for many different quivers (see section 3 of [GLS05] for more details). Now theorem 3 tells us that if $\Lambda_n$ and $\Lambda_m$ are separably equivalent then the preprojective algebras of type $A_{n-1}$ and $A_{m-1}$ are separably equivalent.

**Theorem 8:**

> Let $\Lambda_n$ denote the truncated polynomial algebra $k[x]/(x^n)$ over an algebraically closed field $k$.
>
> The algebras $\Lambda_n$ and $\Lambda_m$ are not separably equivalent for positive integers $n \leq 6$ and $m \neq n$.

*Proof*[1].   We prove the claim by handling each $n$ on a case-by-case basis.

$n = 1$

The algebra $\Lambda_1$ is isomorphic to the field $k$ and so if $\Lambda_m$ is separably equivalent to $\Lambda_1$ then $\Lambda_m$ must be separable and hence semisimple. It is clear that for any $m > 1$ the algebra $\Lambda_m$ is not semisimple and hence the algebras are separably inequivalent.

[GLS05] Geiss, Leclerc, and Schröer, *Semicanonical bases and preprojective algebras*, Ann. Sci. École Norm. Sup. (4) **38** (2005), no. 2, 193–253

---

[1]Note that Linckelmann gives an alternative proof of the $n = 2$ case in example 10.7 of [Lin11a].

[Lin11a] Linckelmann, *Cohomology of block algebras of finite groups*, Representations of algebras and related topics, EMS Ser. Congr. Rep., Eur. Math. Soc., Zürich, 2011, pp. 189–250

$n = 2$   If $\Lambda_2$ is separably equivalent to $\Lambda_m$ for some $m > 2$ then by the discussion above we see that this gives a separable equivalence of the preprojective algebras of types $A_1$ and $A_{m-1}$. The first of these is isomorphic to the base field $k$, but it is clear that for $m > 2$ the preprojective algebra $A_{m-1}$ is not semisimple, therefore the algebras cannot be separably equivalent.

$n = 3$

We denote by $\Gamma_n$ the preprojective algebra of type $A_n$. As in the last two cases $\Lambda_3$ and $\Lambda_m$ being separably equivalent means that $\Gamma_2$ and $\Gamma_{m-1}$ are separably equivalent. From theorem 3 we can see that the two categories $\text{Fun}(\underline{\text{mod}}\,\Gamma_2, \text{mod}\,k)$ and $\text{Fun}(\underline{\text{mod}}\,\Gamma_{m-1}, \text{mod}\,k)$ are separably equivalent. Now $\Gamma_2$ is the path algebra of

$$\bullet \underset{\beta}{\overset{\alpha}{\rightleftarrows}} \bullet \qquad \alpha\beta = \beta\alpha = 0$$

with the projective modules given by

$$k \underset{0}{\overset{1}{\rightleftarrows}} k \qquad\qquad k \underset{1}{\overset{0}{\rightleftarrows}} k$$

Thus in $\underline{\text{mod}}\,\Gamma_2$, once we have factored out projective modules, we are left with just the simple modules

$$k \rightleftarrows 0 \qquad\qquad 0 \rightleftarrows k$$

with no non-trivial maps between them. This shows that $\text{Fun}(\underline{\text{mod}}\,\Gamma_2, \text{mod}\,k)$ is equivalent to the representations of the quiver with two vertices and no arrows: a semisimple algebra.

To show that $\text{Fun}(\underline{\text{mod}}\,\Gamma_{m-1}, \text{mod}\,k)$ is not semisimple when $m > 2$ we need only demonstrate that the AR-quiver of $\Gamma_{m-1}$ contains an arrow between two vertices representing non-projective modules. When $m > 5$ the AR-quiver of $\Gamma_{m-1}$ is infinite (see [DR92, proposition 6.3]) and thus such an arrow must exist. For smaller $m$, section 20 of [GLS05] gives the explicit AR-quivers and shows that the AR-quivers of both $\Gamma_3$ and $\Gamma_4$ contain an arrow between non-projectives.

$n = 4$

If $\Lambda_4$ and $\Lambda_m$ were separably equivalent then their tensor products with a fixed third algebra would be equivalent also. For an algebra $A$, the algebra of upper triangular matrices with entries from $A$

$$\begin{pmatrix} A & A \\ 0 & A \end{pmatrix}$$

[DR92] Dlab and Ringel, *The module theoretical approach to quasi-hereditary algebras*, Representations of algebras and related topics (Kyoto, 1990), London Math. Soc. Lecture Note Ser., vol. 168, Cambridge Univ. Press, Cambridge, 1992, pp. 200–224

[GLS05] Geiss, Leclerc, and Schröer, *Semicanonical bases and preprojective algebras*, Ann. Sci. École Norm. Sup. (4) **38** (2005), no. 2, 193–253

is called the 2 × 2 triangular matrix algebra, $T_2(A)$. This is isomorphic to the tensor product of $A$ with the path algebra 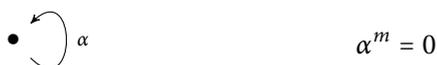 via the isomorphism

$$
\begin{array}{rcl}
A \otimes \left(1 \xrightarrow{\alpha} 2\right) & \longrightarrow & T_2(A) \\
a \otimes e_1 & \mapsto & \begin{pmatrix} a & 0 \\ 0 & 0 \end{pmatrix} \\
a \otimes e_2 & \mapsto & \begin{pmatrix} 0 & 0 \\ 0 & a \end{pmatrix} \\
a \otimes \alpha & \mapsto & \begin{pmatrix} 0 & a \\ 0 & 0 \end{pmatrix}
\end{array}
$$

The representation types of algebras of this form were classified in [LS00]. These were classified via lists of quivers that may appear as a factor algebra of a subquiver of a Galois covering. The relevant sections are theorems 1 and 4, together with the lists of quivers in sections 2 and 5. The algebra $\Lambda_m$ is the path algebra for the quiver

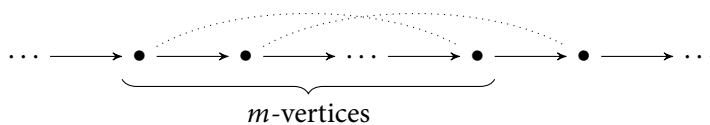

$\alpha^m = 0$

[LS00] Leszczyński and Skowroński, *Tame triangular matrix algebras*, Colloq. Math. **86** (2000), no. 2, 259–303

The Galois covering of this is given by the quiver

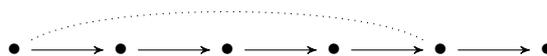

$m$-vertices

with dotted lines representing the relation that the path is zero. For more details on this example and Galois coverings in general see [Gab81, 2.8ff.]. The key point to note regarding this quiver is that for $m > 4$ the Galois covering contains a subquiver with

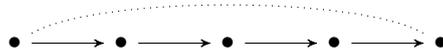

[Gab81] Gabriel, *The universal cover of a representation-finite algebra*, Representations of algebras (Puebla, 1980), Lecture Notes in Math., vol. 903, Springer, Berlin-New York, 1981, pp. 68–105

as a factor algebra, which is [LS00, 2.74] and this means that $T_2(\Lambda_m)$ is wild. When $m = 4$ there is no subquiver containing a factor algebra of wild type but it does contain

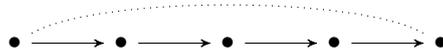

as a factor algebra, which is [LS00, 5.12] and shows that $T_2(\Lambda_4)$ is tame. Finally we can check that no quiver in section 5 of the text appears as a subfactor of the Galois quiver when $m < 4$ and hence for these $T_2(\Lambda_m)$ are of finite type. We have

$T_2(\Lambda_n)$ has    finite representation type for    $n < 4$  
               tame                                             $n = 4$  
               wild                                              $n > 4$

Since $T_2(\Lambda_4)$ and $T_2(\Lambda_m)$ have different representation type for $m > 4$ theorem 6 tells us that $\Lambda_4$ and $\Lambda_m$ cannot be separably equivalent. Note that this method also gives an alternative proof for the case $n = 3$.

$n = 5$ and $n = 6$

Finally we can see from [DR92, proposition 6.3] and [GLS05, proposition 3.3] that

$$\Gamma_n \text{ has } \begin{array}{ll} \text{finite representation type for} & n < 5 \\ \text{tame} & n = 5 \\ \text{wild} & n > 5 \end{array}$$

and since a separable equivalence of $\Lambda_n$ and $\Lambda_m$ induces an equivalence of $\Gamma_{n-1}$ and $\Gamma_{m-1}$ we see that when $n \in \{5, 6\}$ and $m \neq n$ then $\Lambda_n$ and $\Lambda_m$ cannot be separably equivalent. □

*Remark.* The assumption of algebraic closure is only required for the arguments regarding representation type, as such for $n \leq 3$ we do not actually require the field to be algebraically closed.

Theorem 8 is a long way from answering the question as to whether or not algebras for cyclic groups can be separably equivalent. It does however demonstrate many examples of how one can show that algebras are not separably equivalent, using many of the propositions of the preceding sections. The proof of the theorem uses representation type to differentiate between algebras for $n \leq 6$. For larger $n$ all the algebras we have constructed from $\Gamma_n$ have wild representation type and so it would appear new methods will be required to show the inequivalence of these algebras. We conclude with the corollaries:

**Corollary.** *Let k be a field of characteristic* 2. *The group algebras $kC_2$, $kC_4$ and $kC_{2^n}$ are pairwise separably inequivalent for any $n > 2$.*

**Corollary.** *Let k be a field of characteristic* 3. *The group algebras $kC_3$ and $kC_{3^n}$ are separably inequivalent for any $n > 1$.*

**Corollary.** *Let k be a field of characteristic* 5. *The group algebras $kC_5$ and $kC_{5^n}$ are separably inequivalent for any $n > 1$.*

# References

[Alp86]   J. L. Alperin, *Local representation theory*, Cambridge Studies in Advanced Mathematics, vol. 11, Cambridge University Press, Cambridge, 1986.

[AE81]   J. L. Alperin and L. Evens, *Representations, resolutions and Quillen's dimension theorem*, J. Pure Appl. Algebra **22** (1981), no. 1, 1–9.

[DR92]   Dlab and Ringel, *The module theoretical approach to quasi-hereditary algebras*, Representations of algebras and related topics (Kyoto, 1990), London Math. Soc. Lecture Note Ser., vol. 168, Cambridge Univ. Press, Cambridge, 1992, pp. 200–224

[GLS05]   Geiss, Leclerc, and Schröer, *Semicanonical bases and preprojective algebras*, Ann. Sci. École Norm. Sup. (4) **38** (2005), no. 2, 193–253

[BE11] Petter Andreas Bergh and Karin Erdmann, *The representation dimension of Hecke algebras and symmetric groups*, Adv. Math. **228** (2011), no. 4, 2503–2521.

[CB91] William Crawley-Boevey, *Tame algebras and generic modules*, Proc. London Math. Soc. (3) **63** (1991), no. 2, 241–265.

[CB92] \_\_\_\_\_\_, *Modules of finite length over their endomorphism rings*, Representations of algebras and related topics (Kyoto, 1990), London Math. Soc. Lecture Note Ser., vol. 168, Cambridge Univ. Press, Cambridge, 1992, pp. 127–184.

[DI71] Frank DeMeyer and Edward Ingraham, *Separable algebras over commutative rings*, Lecture Notes in Mathematics, vol. 181, Springer-Verlag, Berlin-New York, 1971.

[DR92] Vlastimil Dlab and Claus Michael Ringel, *The module theoretical approach to quasi-hereditary algebras*, Representations of algebras and related topics (Kyoto, 1990), London Math. Soc. Lecture Note Ser., vol. 168, Cambridge Univ. Press, Cambridge, 1992, pp. 200–224.

[Dro77] Ju. A. Drozd, *Tame and wild matrix problems*, Matrix problems (Russian), Akad. Nauk Ukrain. SSR Inst. Mat., Kiev, 1977, pp. 104–114.

[Dro80] \_\_\_\_\_\_, *Tame and wild matrix problems*, Representation theory, II (Proc. Second Internat. Conf., Carleton Univ., Ottawa, Ont., 1979), Lecture Notes in Math., vol. 832, Springer, Berlin, 1980, pp. 242–258.

[Gab81] P. Gabriel, *The universal cover of a representation-finite algebra*, Representations of algebras (Puebla, 1980), Lecture Notes in Math., vol. 903, Springer, Berlin-New York, 1981, pp. 68–105.

[GLS05] Christof Geiss, Bernard Leclerc, and Jan Schröer, *Semicanonical bases and preprojective algebras*, Ann. Sci. École Norm. Sup. (4) **38** (2005), no. 2, 193–253.

[LS00] Zbigniew Leszczyński and Andrzej Skowroński, *Tame triangular matrix algebras*, Colloq. Math. **86** (2000), no. 2, 259–303.

[Lin11a] Markus Linckelmann, *Cohomology of block algebras of finite groups*, Representations of algebras and related topics, EMS Ser. Congr. Rep., Eur. Math. Soc., Zürich, 2011, pp. 189–250.

[Lin11b] \_\_\_\_\_\_, *Finite generation of Hochschild cohomology of Hecke algebras of finite classical type in characteristic zero*, Bull. Lond. Math. Soc. **43** (2011), no. 5, 871–885.

[Pre09] Mike Prest, *Purity, spectra and localisation*, Encyclopedia of Mathematics and its Applications, vol. 121, Cambridge University Press, Cambridge, 2009.

School of Mathematics and The Heilbronn Institute for Mathematical Research, University of Bristol, Bristol, BS8 1TW
  *E-mail address*: simon.peacock@bristol.ac.uk